\pdfoutput=1
\documentclass[11pt]{article}
\usepackage[left=1in,right=1in,top=1in,bottom=1in]{geometry}
\usepackage{times}
\usepackage{expl3}
\usepackage{cite}
\usepackage[table]{xcolor}
\usepackage{multirow}
\usepackage{stackengine} 
\usepackage{hhline}
\usepackage{lipsum}
\usepackage{titlesec}
\usepackage{wrapfig}
\usepackage{enumerate}
\usepackage{epsfig}
\usepackage{amsmath}
\usepackage{tabularx}
\usepackage{array}
\usepackage{booktabs}
\usepackage{enumitem}
\usepackage{bbm}
\usepackage{calc}
\usepackage{graphicx}
\usepackage{amsmath}
\usepackage[title]{appendix}
\usepackage{amssymb}
\usepackage{epstopdf}
\usepackage{boldline}
\usepackage{arydshln}
\usepackage{calligra}
\usepackage{bm}
\usepackage{url}
\usepackage{blindtext}
\usepackage{accents}

\newcommand{\define}{\stackrel{\mbox{\tiny def}}{=}}

\newtheorem{definition}{Definition}
\newtheorem{theorem}{Theorem}

\newtheorem{lemma}{Lemma}

\newtheorem{remark}{Remark}
\usepackage{mathtools}
\usepackage{epstopdf}
\usepackage{balance}
\usepackage{thmtools}
\usepackage{thm-restate}
\usepackage{hyperref}
\usepackage{cleveref}
\usepackage[mathscr]{euscript}

\usepackage[ruled,vlined]{algorithm2e}
\include{pythonlisting}
\newcommand{\ostar}{\mathbin{\mathpalette\make@circled\star}}

\makeatletter
\newcommand{\removelatexerror}{\let\@latex@error\@gobble}
\makeatother
\makeatletter
\newcommand*{\rom}[1]{\expandafter\@slowromancap\romannumeral #1@}
\makeatother

\ExplSyntaxOn
\newcommand\latinabbrev[1]{
  \peek_meaning:NTF . {
    #1\@}%
  { \peek_catcode:NTF a {
      #1.\@ }%
    {#1.\@}}}
\ExplSyntaxOff

\titleclass{\subsubsubsection}{straight}[\subsubsection]

\begin{document}
\vspace{1cm}
\title{Tail Bounds for Functions of Weighted Tensor Sums Derived from Random Walks on Riemannian Manifolds}\vspace{1.8cm}
\author{Shih-Yu~Chang
\thanks{Shih-Yu Chang is with the Department of Applied Data Science,
San Jose State University, San Jose, CA, U. S. A. (e-mail: {\tt
shihyu.chang@sjsu.edu}). 
           }}

\maketitle

\begin{abstract}
This paper presents significant advancements in tensor analysis and the study of random walks on manifolds. It introduces new tensor inequalities derived using the Mond-Pecaric method, which enriches the existing mathematical tools for tensor analysis. This method, developed by mathematicians Mond and Pecaric, is a powerful technique for establishing inequalities in linear operators and matrices, using functional analysis and operator theory principles. The paper also proposes novel lower and upper bounds for estimating column sums of transition matrices based on their spectral information, which is critical for understanding random walk behavior. Additionally, it derives bounds for the right tail of weighted tensor sums derived from random walks on manifolds, utilizing the spectrum of the Laplace-Beltrami operator over the underlying manifolds and new tensor inequalities to enhance the understanding of these complex mathematical structures.
\end{abstract}


\section{Introduction}\label{sec: Introduction}

Recently, the study about random tensors has attracted significant academic interest due to their critical role in advancing our understanding of multi-dimensional data structures and their applications across different science and engineering fields. Tensors generalize matrices to higher dimensions, making them ideal for representing complex data sets that arise in fields like machine learning, signal processing, quantum mechanics, and neuroscience. People can understand the probabilistic behavior of these higher-dimensional objects, gaining insights into their typical properties and performance in multi-dimensional data analysis via studying random tensors~\cite{carrozza2016n}. Random tensors also enable the development of theoretical frameworks for understanding tensor decompositions, spectral properties, and optimization in higher dimensions, which are essential for tackling problems in areas such as tensor completion, multi-way data analysis, and the study of complex networks. Moreover, the study of random tensors helps in addressing challenges related to computational complexity, providing bounds and concentration inequalities that are crucial for ensuring the stability and reliability of tensor-based algorithms. Ultimately, random tensors serve as a foundational tool in both theoretical and applied mathematics, offering deep insights into the behavior of high-dimensional systems and contributing to the development of robust methodologies for handling complex data~\cite{ouerfelli2022random}.

Considering the tail bounds of the sum of random variables is an old problem in probability theory. The original purpose of studying tail bunds of random variablesis to understand the likelihood of extreme deviations from expected outcomes. Tail bounds provide probabilistic guarantees that the sum will remain within a certain range, which is crucial in fields like finance, machine learning, and statistical analysis where risk management and error estimation are vital. These bounds help ensure the stability and reliability of algorithms, particularly in high-dimensional settings, and are key to deriving concentration inequalities that underpin the performance and safety of complex systems involving randomness~\cite{tropp2015introduction}. The study in tail bounds also brings great impacts in various areas in mathematics. In~\cite{ledoux2001concentration}. The author explores the fundamental techniques and examples of the concentration of measure phenomenon, which was introduced in the early 1970s by V. Milman in the context of the asymptotic geometry of Banach spaces. The concept about concentration over the metric-measure space has gained significant interest across various fields, including geometry, discrete mathematics, complexity theory, functional analysis, infinite-dimensional integration, and probability theory itself. The extension of tail bounds from a sequence of random variables to a sequence of random tensors based on indepedent assumption is made by the following works~\cite{chang2022convenient,chang2022generalizedHW,chang2022randomDTI,chang2022randomDTI_App,chang2022randomMOI,chang2023tail,chang2023randomB_I,chang2023randomB_II,chang2023tail_RP}.

Not all sequences of random objects are independent because real-world processes often exhibit dependencies and varying distributions. For example, in time series data, each observation may depend on previous ones, creating correlations that violate independence. Non independent sequences of random objects are common in applications like finance, where market conditions fluctuate, or in sensor networks, where measurements might be influenced by environmental changes. Therefore, it is important to consider non independent sequences of random objects for accurately modeling and analyzing complex, dynamic systems. Following this spirit, we extend previous works about tail bounds from independent random tensors to non-independent random tensors by utilizing the random walk model ~\cite{chang2021TensorExp,chang2021T-TensorExp}. Those tail bounds provided by~\cite{chang2021TensorExp,chang2021T-TensorExp} are controlled by the spectrum of the Laplacian matrix of the graph, which is the underlying space of the random walk model. 

In~\cite{chang2024chernoff_I}, we first attempt to study non-independent random tensors dervied from random walks over over manifolds. There are three main reasons to study specturm of Laplace-Beltrami operator over a manifold. First, the eigenvalues of the Laplace-Beltrami operator are intrinsically linked to the geometry of the manifold. For instance, they reflect the manifold's volume, curvature, and the way it stretches or compresses. This makes the spectrum a powerful tool for distinguishing between different manifolds, even when they share some superficial similarities. Second, the spectrum is connected to the manifold's topology. Certain topological features, such as the number of holes or connected components, can influence the spectrum, thereby allowing mathematicians to infer topological information from the eigenvalues. Third, the Laplace-Beltrami spectrum has applications in various fields beyond pure mathematics, including physics, where it relates to the behavior of waves and heat distribution on surfaces, and computer science, particularly in areas like shape analysis and machine learning. By studying this spectrum, researchers can better understand complex manifolds, solve differential equations on them, and apply these insights across diverse scientific and engineering disciplines. Therefore, the tail bounds given by~\cite{chang2024chernoff_I} are associated with spectrum information of the underlying manifold. Our strategy is to approximate the underlying manifold by a graph that has similar spectrum information to the spectrum of the Laplace-Beltrami operator based on the work from~\cite{burago2015graph}. Then, we can derive the tensor Chernoff bound and establish its range for random walks on a Riemannian manifold according to the underlying manifold’s spectral characteristics.

This paper makes significant contributions to tensor analysis and random walks on manifolds. It first introduces new tensor inequalities using the Mond-Pecaric method~\cite{pecaric2005mond,fujii2012recent}, enriching the mathematical tools available for tensor analysis. The Mond-Pecaric method, named after mathematicians Mond and Pecaric, is a powerful technique for deriving inequalities, particularly in linear operators and matrices. This approach constructs operator inequalities through functional analysis and operator theory, utilizing properties like convexity and monotonicity. Widely applied in matrix and tensor analysis, it helps establish precise and general inequalities relevant to optimization and numerical analysis. The paper also presents novel lower and upper bounds for estimating column sums of transition matrices based on their spectral information, crucial for understanding random walk behavior through their transition matrices. Finally, it derives lower and upper bounds for the right tail of weighted tensor sums derived from random walks on manifolds, leveraging new tensor inequalities and random walks transition matrix spectrum estimation using the Laplace-Beltrami operator's spectrum on the underlying manifold.

The rest of this paper is organized as follows. In Section~\ref{sec: Tensor Inequalities by Mond-Pecaric Method}, we will provide new tensor inequalities according to the Mond-Pecaric method. In Section~\ref{sec: Column Sums Estimation by Spectrum of Transition Matrix}, we will derive lower bounds and upper bounds estimation for column sums of the transition matrix in terms of the spectrum information of the transition matrix. The transition matrix and its spectrum information for random walks over manifolds is explored in Section~\ref{sec: Random Walk Transition Matrix Spectrum Approximation by Laplace-Beltrami Operator Spectrum}. Finally, in Section~\ref{sec: Tail Bounds for Functions of Weighted Tensor Sums}, we will derive tail bounds for functions of weighted tensor sums derived from random walks on Riemannian manifolds. 

\noindent \textbf{Nomenclature:} To simplify notation, let $\mathbb{I}_1^M$ be defined as $\prod_{i=1}^M I_i$, where $I_i$ indicates the size of the $i$-th dimension of a tensor. The notions about Hermitian tensor and its eigenvalues and eigen-tensor is given by~\cite{chang2022convenient}. $\lambda_{i}(\mathcal{X})$ represents the $i$-th eigenvalue of the tensor $\mathcal{X}$, where the tensor $\mathcal{X}$ becomes a matrix if the order of the tensor $\mathcal{X}$ is 2. $\geq$ and $\leq$ between tensors are Loewner orderings. We also assume that all tensor-valued functions of Hermitian tensors are Hermitian tensors in this work. The script $\mathrm{H}$ represents the Hermitian operation. $\left\Vert\mathcal{H}\right\Vert$ represents any unitarily invariant tensor norms~\cite{chang2021TensorExp}. 

\section{Tensor Inequalities by Mond-Pecaric Method}\label{sec: Tensor Inequalities by Mond-Pecaric Method}

The purpose of this section is to establish Theorem~\ref{thm: Tensor Inequalities by MP}, which will provide new tensor inequalities obtained by applying the Mond-Pecaric method. 

We will begin by providing Definition~\ref{def: nomralized positive linear map for tensors} about a \emph{normalized positive linear map for tensors}, which are the finite dimensional operators, a special case of \emph{normalized positive linear map for operators} given by~\cite{pecaric2005mond,fujii2012recent}.
\begin{definition}\label{def: nomralized positive linear map for tensors}
Let $\mathscr{A}(\mathfrak{H})$ and $\mathscr{A}(\mathfrak{K})$ be semi-algebras of all tensors on the vector space $\mathfrak{H}$ and on the vector space $\mathfrak{K}$, respectively. The symbols $\mathcal{I}_{\mathfrak{H}}$ and $\mathcal{I}_{\mathfrak{K}}$ are represented as identity tensors for the vector space $\mathfrak{H}$ and the vector space $\mathfrak{K}$, respectively. A normalized positive linear map for tensors is a map $\Psi: \mathscr{A}(\mathfrak{H}) \rightarrow \mathscr{A}(\mathfrak{K})$ such that the following three conditions are satisfied:
\begin{enumerate}
\item \textbf{Linear:} we have 
\begin{eqnarray}\label{eq: linerity}
\Psi(a\mathcal{X}+b\mathcal{Y})=a\Psi(\mathcal{X})+b\Psi(\mathcal{Y}), 
\end{eqnarray}
for any $a,b \in \mathbb{C}$ and $\mathcal{X},\mathcal{Y} \in \mathscr{A}(\mathfrak{H})$. For a map $\Psi$ that satisies Eq.~\eqref{eq: linerity} is named as a \emph{linear} map. 
\item \textbf{Normalized:} For any linear map $\Psi$ that satisfies $\Psi(\mathcal{I}_{\mathfrak{H}}) \rightarrow \mathcal{I}_{\mathfrak{K}}$ is calles as a \emph{normalized} map.   
\item \textbf{Positive:} For any linear map $\Psi$ that satisfies $\Psi(\mathcal{X})\geq\Psi(\mathcal{Y})$ for any $\mathcal{X}\geq\mathcal{Y}$ is calles as a \emph{positive} map.   
\end{enumerate}
\end{definition}

There are many possible ways to construct the map $\Psi$ that satisfies Definition~\ref{def: nomralized positive linear map for tensors}. In this work, we will adopt those $\Psi$ maps that keeps Hermitian property, i.e., $ \Psi(\mathcal{X})$ will be Hermitian tensor if $\mathcal{X}$ is a Hermitian tensor. For example, we can select $\Psi(\mathcal{X}) = \mathcal{U}\star\mathcal{X}\star\mathcal{U}^{\mathrm{H}}$, where $\mathcal{U}\star\mathcal{U}^{\mathrm{H}}=\mathcal{I}$.  

Lemma~\ref{lma: func bounded by lines} is provided to bound a convex function.
\begin{lemma}\label{lma: func bounded by lines}
Given a convex function $g$ in the real interval $[c, d]$, we have 
\begin{eqnarray}\label{eqU: lma: func bounded by lines}
g(s)&\leq&  \frac{g(d)-g(c)}{d-c}s + \frac{dg(c)-cg(d)}{d-c},
\end{eqnarray}
and
\begin{eqnarray}\label{eqL: lma: func bounded by lines}
g(s)&\geq&\frac{g(d)-g(c)}{d-c}s + \left[g\left((g')^{-1}\left(\frac{g(d)-g(c)}{d-c}\right)\right)- \frac{g(d)-g(c)}{d-c}(g')^{-1}\left(\frac{g(d)-g(c)}{d-c}\right)\right],
\end{eqnarray}
where $(g')^{-1}$ is the inverse function with respect to the first derivative of the function $g$. 
\end{lemma}
\textbf{Proof:}
Because the function $g$ is a convex function in the real interval $[c, d]$, we can upper bound this function $g$ by a linear function passing points $(c,g(c))$ and $(d,g(d))$. The equation for this line that passes points $(c,g(c))$ and $(d,g(d))$ can be expressed by
\begin{eqnarray}\label{eqU1: lma: func bounded by lines}
\frac{g(d)-g(c)}{d-c}s + \frac{dg(c)-cg(d)}{d-c}.
\end{eqnarray}
Then, we have the inequality provided by Eq.~\eqref{eqU: lma: func bounded by lines}. 

On the other hand, it is important to note that the value of the function $g$ within the interval $[c, d]$ will always be greater than or equal to the value of a linear function with a slope of $\frac{g(d) - g(c)}{d - c}$ that passes through the point $s_0$ with $g'(s_0) = \frac{g(d) - g(c)}{d - c}$. The equation for this line can be expressed by
\begin{eqnarray}\label{eqL1: lma: func bounded by lines}
\frac{g(d)-g(c)}{d-c}s + \left[g\left((g')^{-1}\left(\frac{g(d)-g(c)}{d-c}\right)\right)- \frac{g(d)-g(c)}{d-c}(g')^{-1}\left(\frac{g(d)-g(c)}{d-c}\right)\right].
\end{eqnarray}
Then, we also have the inequality given by Eq.~\eqref{eqL: lma: func bounded by lines}. 
$\hfill \Box$

For notation simplicity, we will have the following abbreviation definitions
\begin{eqnarray}\label{eq: m g b g U b g L def}
m_g &\define& \frac{g(d)-g(c)}{d-c}, \nonumber \\
b_{g,U}&\define&\frac{dg(c)-cg(d)}{d-c}, \nonumber \\
b_{g,L}&\define& \left[g\left((g')^{-1}\left(\frac{g(d)-g(c)}{d-c}\right)\right)- \frac{g(d)-g(c)}{d-c}(g')^{-1}\left(\frac{g(d)-g(c)}{d-c}\right)\right]. 
\end{eqnarray}

Following Lemma~\ref{lma: ext Thm 2.4} is the extension of Theorem 2.4 in~\cite{pecaric2005mond} which will be used in proving Theorem~\ref{thm: Tensor Inequalities by MP}.

\begin{lemma}\label{lma: ext Thm 2.4}
We are given $\ell$ Hermitian tensors $\mathcal{X}_i$ with eigenvalues within in the real interval $[c, d]$, $\ell$  normalized positive linear maps for tensors $\Psi_i$, and a probability vector with entries $w_i$, i.e., $\sum\limits_{i=1}^\ell w_i = 1$ and $w_i \geq 0$, for $i=1,2,\ldots ,\ell$. For continuous real functions $g,h$  with $g$ convexity in the real interval $[c, d]$ and any real number $c_r$, we have
\begin{eqnarray}\label{eqU1:lma: ext Thm 2.4}
\sum\limits_{i=1}^{\ell}w_i \Psi_i (g(\mathcal{X}_i))&\leq&c_r h\left(\sum\limits_{i=1}^\ell w_i \Psi_i(\mathcal{X}_i)\right)+ \max\limits_{c \leq s \leq d}\left[ m_g s + b_{g,U} - c_r h(s)\right]\mathcal{I}_{\mathfrak{K}}.
\end{eqnarray}
On the other hand, we also have 
\begin{eqnarray}\label{eqL1:lma: ext Thm 2.4}
\sum\limits_{i=1}^{\ell}w_i \Psi_i (g(\mathcal{X}_i))&\geq&c_r h\left(\sum\limits_{i=1}^\ell w_i \Psi_i(\mathcal{X}_i)\right)+ \min\limits_{c \leq s \leq d}\left[m_g s +b_{g,L}- c_r h(s)\right]\mathcal{I}_{\mathfrak{K}}.
\end{eqnarray}
\end{lemma}
\textbf{Proof:}
From Eq.~\eqref{eqU: lma: func bounded by lines} in Lemma~\ref{lma: func bounded by lines} and spectral mapping theorem, we have
\begin{eqnarray}\label{eqU2:lma: ext Thm 2.4}
g(\mathcal{X}_i)&\leq& m_g \mathcal{X}_i + b_{g,U} \mathcal{I}_{\mathfrak{H}},
\end{eqnarray}
where $i=1,2,\ldots,\ell$. By apply the map $\Psi_i$ to both sides of Eq.~\eqref{eqU2:lma: ext Thm 2.4}, we have
\begin{eqnarray}\label{eqU3:lma: ext Thm 2.4}
\Psi_i(g(\mathcal{X}_i))&\leq& \Psi_i(m_g \mathcal{X}_i + b_{g,U} \mathcal{I}_{\mathfrak{H}}) \nonumber \\
&=& m_g \Psi_i(\mathcal{X}_i) + b_{g,U}\Psi_i (\mathcal{I}_{\mathfrak{H}}) \nonumber \\
&=& m_g \Psi_i(\mathcal{X}_i) + b_{g,U}\mathcal{I}_{\mathfrak{K}} 
\end{eqnarray}
where $i=1,2,\ldots,\ell$. If we multiply the weight $w_i$ to both sides of Eq.~\eqref{eqU3:lma: ext Thm 2.4} and sum over the index $i$, we obtain:
\begin{eqnarray}\label{eqU4:lma: ext Thm 2.4}
\sum\limits_{i=1}^{\ell}w_i\Psi_i(g(\mathcal{X}_i))&\leq& m_g\sum\limits_{i=1}^{\ell}w_i\Psi_i(\mathcal{X}_i) + b_{g,U}\mathcal{I}_{\mathfrak{K}}.
\end{eqnarray}
For any given real number $c_r$ and Eq.~\eqref{eqU4:lma: ext Thm 2.4}, we have 
\begin{eqnarray}\label{eqU5:lma: ext Thm 2.4}
\sum\limits_{i=1}^{\ell}w_i\Psi_i(g(\mathcal{X}_i)) - c_r h\left(\sum\limits_{i=1}^{\ell}w_i\Psi_i(\mathcal{X}_i)\right)&\leq& m_g\sum\limits_{i=1}^{\ell}w_i\Psi_i(\mathcal{X}_i) + b_{g,U}\mathcal{I}_{\mathfrak{K}}  - c_r h\left(\sum\limits_{i=1}^{\ell}w_i\Psi_i(\mathcal{X}_i)\right)\nonumber \\
&\leq_1& \left\{\max\limits_{c \leq s \leq d}\left[m_g s +b_{g,U} - c_r h(s)\right]\right\}\mathcal{I}_{\mathfrak{K}}.
\end{eqnarray}
where $\leq_1$ comes from the fact that the range for eigenvalues of  $\sum\limits_{i=1}^{\ell}w_i\Psi_i(\mathcal{X}_i)$ is also with in the real interval $[c,d]$. Then, we have the inequality given by Eq.~\eqref{eqU1:lma: ext Thm 2.4} by rearranging the term in Eq.~\eqref{eqU5:lma: ext Thm 2.4}.

From Eq.~\eqref{eqL: lma: func bounded by lines} in Lemma~\ref{lma: func bounded by lines} and spectral mapping theorem, we have
\begin{eqnarray}\label{eqL2:lma: ext Thm 2.4}
g(\mathcal{X}_i)&\geq& m_g \mathcal{X}_i + b_{g,L} \mathcal{I}_{\mathfrak{H}},
\end{eqnarray}
where $i=1,2,\ldots,\ell$. By appling the map $\Psi_i$ to both sides of Eq.~\eqref{eqL2:lma: ext Thm 2.4}, we have
\begin{eqnarray}\label{eqL3:lma: ext Thm 2.4}
\Psi_i(g(\mathcal{X}_i))&\geq& \Psi_i(m_g \mathcal{X}_i + b_{g,L} \mathcal{I}_{\mathfrak{H}}) \nonumber \\
&=& m_g \Psi_i(\mathcal{X}_i) + b_{g,L}\Psi_i (\mathcal{I}_{\mathfrak{H}}) \nonumber \\
&=& m_g \Psi_i(\mathcal{X}_i) + b_{g,L}\mathcal{I}_{\mathfrak{K}} 
\end{eqnarray}
where $i=1,2,\ldots,\ell$. If we multiply the weight $w_i$ to both sides of Eq.~\eqref{eqL3:lma: ext Thm 2.4} and sum over the index $i$, we obtain:
\begin{eqnarray}\label{eqL4:lma: ext Thm 2.4}
\sum\limits_{i=1}^{\ell}w_i\Psi_i(g(\mathcal{X}_i))&\geq& m_g\sum\limits_{i=1}^{\ell}w_i\Psi_i(\mathcal{X}_i) + b_{g,L}\mathcal{I}_{\mathfrak{K}}.
\end{eqnarray}
For any given real number $c_r$ and Eq.~\eqref{eqU4:lma: ext Thm 2.4}, we have 
\begin{eqnarray}\label{eqL5:lma: ext Thm 2.4}
\sum\limits_{i=1}^{\ell}w_i\Psi_i(g(\mathcal{X}_i)) - c_r h\left(\sum\limits_{i=1}^{\ell}w_i\Psi_i(\mathcal{X}_i)\right)&\geq& m_g\sum\limits_{i=1}^{\ell}w_i\Psi_i(\mathcal{X}_i) + b_{g,L}\mathcal{I}_{\mathfrak{K}}  - c_r h\left(\sum\limits_{i=1}^{\ell}w_i\Psi_i(\mathcal{X}_i)\right)\nonumber \\
&\geq_1& \left\{\min\limits_{c \leq s \leq d}\left[m_g s +b_{g,L} - c_r h(s)\right]\right\}\mathcal{I}_{\mathfrak{K}}.
\end{eqnarray}
where $\geq_1$ comes from the fact that the range for eigenvalues of  $\sum\limits_{i=1}^{\ell}w_i\Psi_i(\mathcal{X}_i)$ is with in the real interval $[c,d]$. Then, we have the inequality given by Eq.~\eqref{eqL1:lma: ext Thm 2.4} by rearranging the term in Eq.~\eqref{eqL5:lma: ext Thm 2.4}.
$\hfill\Box$

We are ready to present Theorem~\ref{thm: Tensor Inequalities by MP} in this section about the upper and the lower bounds for the function of ensemble tensors.

\begin{theorem}\label{thm: Tensor Inequalities by MP}
We are given $\ell$ Hermitian tensors $\mathcal{X}_i$ with eigenvalues within in the real interval $[c, d]$, $\ell$  normalized positive linear maps for tensors $\Psi_i$, and a probability vector with entries $w_i$, i.e., $\sum\limits_{i=1}^\ell w_i = 1$ and $w_i \geq 0$, for $i=1,2,\ldots ,\ell$. For continuous real functions $g,h$  with $g$ convexity in the real interval $[c, d]$ and assume $h(s)>0, m_g s + b_{g,U}>0, m_g s + b_{g,L}>0$ for $s \in [c, d]$, we have
\begin{eqnarray}\label{eq1: thm: Tensor Inequalities by MP}
\frac{1}{\max\limits_{c \leq s \leq d}\left[\frac{m_g s + b_{g,U}}{h(s)}\right]}\sum\limits_{i=1}^{\ell}w_i \Psi_i (g(\mathcal{X}_i))&\leq&h\left(\sum\limits_{i=1}^\ell w_i \Psi_i(\mathcal{X}_i)\right)\nonumber \\
&\leq& \frac{1}{\min\limits_{c \leq s \leq d}\left[\frac{m_g s + b_{g,L}}{h(s)}\right]}\sum\limits_{i=1}^{\ell}w_i \Psi_i (g(\mathcal{X}_i)).
\end{eqnarray}
\end{theorem}
\textbf{Proof:}
We will prove the lower bound for Eq.~\eqref{eq1: thm: Tensor Inequalities by MP} first. From Eq.~\eqref{eqU1:lma: ext Thm 2.4} in Lemma~\ref{lma: ext Thm 2.4}, we have 
\begin{eqnarray}\label{eq2: thm: Tensor Inequalities by MP}
\sum\limits_{i=1}^{\ell}w_i \Psi_i (g(\mathcal{X}_i))&\leq&c_r h\left(\sum\limits_{i=1}^\ell w_i \Psi_i(\mathcal{X}_i)\right)+ \max\limits_{c \leq s \leq d}\left[ m_g s + b_{g,U} - c_r h(s)\right]\mathcal{I}_{\mathfrak{K}}.
\end{eqnarray}
The lower bound is obtained by selecting $c_r = \frac{m_g s + b_{g,U}}{h(s)}$. 

On the other hand, we will prove the upper bound for Eq.~\eqref{eq1: thm: Tensor Inequalities by MP}.  From Eq.~\eqref{eqL1:lma: ext Thm 2.4} in Lemma~\ref{lma: ext Thm 2.4}, we have 
\begin{eqnarray}\label{eq3: thm: Tensor Inequalities by MP}
\sum\limits_{i=1}^{\ell}w_i \Psi_i (g(\mathcal{X}_i))&\geq&c_r h\left(\sum\limits_{i=1}^\ell w_i \Psi_i(\mathcal{X}_i)\right)+ \min\limits_{c \leq s \leq d}\left[ m_g s + b_{g,L} - c_r h(s)\right]\mathcal{I}_{\mathfrak{K}}.
\end{eqnarray}
The upper bound is obtained by selecting $c_r = \frac{m_g s + b_{g,L}}{h(s)}$.
$\hfill \Box$

\section{Column Sums Estimation by Spectrum of Transition Matrix}\label{sec: Column Sums Estimation by Spectrum of Transition Matrix}

In this section, we will derive lower bounds and upper bounds estimation for column sums of the transition matrix in terms of the spectrum information of the transition matrix. We will present two lemmas for the upper bounds estimation and two lemmas for the lower bounds estimation based on different transition matrix properties. 

Given a Markov transition matrix $\bm{P}=[p_{i,j}] \in \mathbb{R}^{N \times N}$ with $p_{i,i} \neq 1$, we define the $j$-th column sum, denoted by $\mathscr{C}_j$, which is $\mathscr{C}_j \define  \sum\limits_{\ell=1}^N p_{\ell,j}$. We also define the $i$-th row sum, denoted by $\mathscr{R}_i$, which is $\mathscr{R}_i \define  \sum\limits_{\ell=1}^N p_{i,\ell}$. Note that $\mathscr{R}_i =1$ since $\bm{P}=[p_{i,j}]$ is a Markov transition matrix. The following two Lemmas will provide the upper bounds for $\mathscr{C}_j$.

\begin{lemma}\label{lma: column sum UB-1}
Given a Markov transition matrix $\bm{P}=[p_{i,j}] \in \mathbb{R}^{N \times N}$ and we assume that $\max\limits_{j \in \{1,2,\ldots N\}}\mathscr{C}_j \geq N-1$, then, we have 
\begin{eqnarray}\label{eq1:lma: column sum UB-1}
\mathscr{C}_j \leq N - \left\vert \lambda_2 \right\vert,
\end{eqnarray}
where $j=1,2,\ldots,N$ and $\left\vert \lambda_2 \right\vert$ is the second largest eigenvalue of the matrix $\bm{P}$.  
\end{lemma}
\textbf{Proof:}
From Corollar 2.2 in~\cite{kirkland2009subdominant}, we know that 
\begin{eqnarray}\label{eq2:lma: column sum UB-1}
\max\limits_{j \in \{1,2,\ldots N\}}\mathscr{C}_j  \leq N - \left\vert \lambda_2 \right\vert.
\end{eqnarray}
This Lemma is proved as $\mathscr{C}_j \leq \max\limits_{j \in \{1,2,\ldots N\}}\mathscr{C}_j$.
$\hfill \Box$

\begin{lemma}\label{lma: column sum UB-2}
Given a Markov transition matrix $\bm{P}=[p_{i,j}] \in \mathbb{R}^{N \times N}$ and we assume that $2p_{i,i}-1 > 0$ for $i=1,2,\ldots,N$, then, for $q \in (0,1)$ we have 
\begin{eqnarray}\label{eq1:lma: column sum UB-2}
\mathscr{C}_j \leq \left[\frac{\min\limits_{i \in \{1,2,\ldots,N\}}\left\vert\lambda_{i}\right\vert}{(2p_{j,j} - 1)^q}\right]^{\frac{1}{1-q}},
\end{eqnarray}
where $j=1,2,\ldots,N$,
\end{lemma}
\textbf{Proof:}
From Corollary A6(a) in~\cite{garren1968bounds}, we have
\begin{eqnarray}\label{eq2:lma: column sum UB-2}
\min\limits_{i \in \{1,2,\ldots,N\}}\left\vert\lambda_{i}\right\vert
&\geq&\left(p_{j,j}-\sum\limits_{k=1, k \neq j}^N p_{j,k}\right)^q
\left(p_{j,j}+\sum\limits_{\ell=1, \ell \neq i}^N p_{\ell,j}\right)^{1-q}\nonumber \\
&=_1&\left(2p_{j,j}-1\right)^q
\mathscr{C}_j^{1-q}
\end{eqnarray}
where $=_1$ comes from the fact that $\bm{P}$ is a Markov transition matrix. This Lemma is obtained by rearrange the terms in Eq.~\eqref{eq2:lma: column sum UB-2}.
$\hfill \Box$

The following two Lemmas will provide the lower bounds for $\mathscr{C}_j$.

\begin{lemma}\label{lma: column sum LB-1}
Given a Markov transition matrix $\bm{P}=[p_{i,j}] \in \mathbb{R}^{N \times N}$ with $p_{i,i} \neq 1$, then, for any $q \in [0,1]$, we have 
\begin{eqnarray}\label{eq1:lma: column sum LB-1}
p_{j,j}+\frac{\left\vert\lambda_j - p_{j,j}\right\vert^{1/(1-q)}}{(1-p_{j,j})^{q/(1-q)}}\leq \mathscr{C}_j
\end{eqnarray}
where $\lambda_j$ is the $j$-th eigenvalue such that 
\begin{eqnarray}\label{eq2:lma: column sum LB-1}
\left\vert\lambda_j - p_{j,j}\right\vert\leq (1-p_{j,j})^q (\mathscr{C}_j-p_{j,j})^{1-q}.
\end{eqnarray}
\end{lemma}
\textbf{Proof:}
According to Ostrowski theorem~\cite{brauer1957theorems}, we have
\begin{eqnarray}
\lambda(\bm{P}) \in \bigcup_{i=1}^{N} S_i, 
\end{eqnarray}
where the set $S_i$ can be expressed by
\begin{eqnarray}
S_i = \{z| \left\vert z - p_{i,i} \right\vert \leq (\mathscr{R}_i-p_{i,i})^q (\mathscr{C}_i-p_{j,j})^{1-q}\}.  
\end{eqnarray}
Therefore, we can find the eigenvalue $\lambda_j$ such that the Eq.~\eqref{eq2:lma: column sum LB-1} is satisfied. This Lemma is proved by rearranging the terms in Eq.~\eqref{eq2:lma: column sum LB-1} since $\mathscr{C}_j - p_{j,j} >0$
$\hfill \Box$

The following Lemma~\ref{lma: column sum LB-2} can be treated as special case of Lemma~\ref{lma: column sum LB-1} by setting the parameter $q$ as $1/2$. 
\begin{lemma}\label{lma: column sum LB-2}
Given a Markov transition matrix $\bm{P}=[p_{i,j}] \in \mathbb{R}^{N \times N}$ with $p_{i,i} \neq 1$, then, we have 
\begin{eqnarray}\label{eq1:lma: column sum LB-2}
\left\vert\lambda_j\right\vert^2\leq\mathscr{C}_j,
\end{eqnarray}
where $\lambda_j$ is the $j$-th eigenvalue such that 
\begin{eqnarray}\label{eq2:lma: column sum LB-2}
\left\vert\lambda_j - p_{j,j}\right\vert\leq \sqrt{(1-p_{j,j})(\mathscr{C}_j-p_{j,j})}.
\end{eqnarray}
\end{lemma}
\textbf{Proof:}
By setting $q=1/2$ in Ostrowski theorem~\cite{brauer1957theorems} and Eq.~\eqref{eq2:lma: column sum LB-1}, we have
\begin{eqnarray}\label{eq2:lma: column sum LB-2}
\left\vert\lambda_j - p_{j,j}\right\vert\leq \sqrt{(1-p_{j,j})(\mathscr{C}_j-p_{j,j})}.
\end{eqnarray}
This implies 
\begin{eqnarray}\label{eq3:lma: column sum LB-2}
\left\vert\lambda_j\right\vert&\leq& p_{j,j} + \sqrt{(1-p_{j,j})(\mathscr{C}_j-p_{j,j})}\nonumber \\
&\leq&\sqrt{\mathscr{C}_j}.
\end{eqnarray}
$\hfill \Box$

\section{Random Walk Transition Matrix Spectrum Approximation by Laplace-Beltrami Operator Spectrum}\label{sec: Random Walk Transition Matrix Spectrum Approximation by Laplace-Beltrami Operator Spectrum}

According to Section 2 of my recent work~\cite{chang2024chernoff_I,burago2015graph}, the following theorem is established.
\begin{theorem}\label{thm: discrete graph Laplacian appro by Laplace-Beltrami Operator on Manifold}
Given a Riemannian manifold $(\mathscr{M}, g)$ of dimension $n$, which consists of a smooth manifold $\mathscr{M}$ equipped with a Riemannian metric $g$, we will construct a weighted graph, denoted by $\mathscr{G}_{\mathscr{M}}(\epsilon,\mu,\kappa) = (\mathscr{V}, \mathscr{E}, \mathscr{W})$to approximate the manifold $\mathscr{M}$ such that  $\mathscr{V} = \{v_i\}$ is the set of vertices sampled from the manifold $\mathscr{M}$ for $i = 1, 2, \ldots, N$, and $\mathscr{E} = \{e_{i,j} = (v_i, v_j)\}$ is the set of edges, satisfying: 
\begin{itemize}
\item All balls with centers $v_i$ and the radius $\epsilon$, denoted by $B_{\epsilon}(v_i)$, can cover $\mathscr{M}$, i.e., $\mathscr{M} \subset \bigcup\limits_{i=1}^N B_{\epsilon}(v_i)$;   
\item The measure function $\mu$ on the set $\{v_i\}$ will allocate the measure $\mu_i$ to the volume of the space $V_i$, where $V_i$ satisfies $\mathscr{M}=\bigcup\limits_{i=1}^N V_i$ and $V_i \subset B_\epsilon(v_i)$; 
\item The edge $e_{i,j}=(v_i, v_j)$ is established if $d_g(v_i, v_j) < \kappa$, and the weight $\omega_{i,j} \subset \mathscr{W}$ for the edge $e_{i,j}\in \mathscr{E}$ is determined by 
\begin{eqnarray}\label{eq: edge weight def}
\omega_{i,j}\define\frac{2(n+2)\Gamma(1+n/2)}{\pi^{n/2}\kappa^{n+2}}\mu_i \mu_j,
\end{eqnarray}
where $\Gamma$ is the Gamma function.
\end{itemize}
Let $\lambda_{\bm{L}_{\mathscr{G}_{\mathscr{M}}},i}$ be the $i$-th eigenvalue of the Laplacian matrix of the graph $\mathscr{G}_{\mathscr{M}}$, i.e., $\bm{L}_{\mathscr{G}_{\mathscr{M}}} = \bm{D}_{\mathscr{G}_{\mathscr{M}}} - \bm{A}_{\mathscr{G}_{\mathscr{M}}}$, where $\bm{D}_{\mathscr{G}_{\mathscr{M}}} $ and $\bm{D}_{\mathscr{G}_{\mathscr{M}}}$ are degree matrix and adjancy matrix of the weighted graph $\mathscr{G}_{\mathscr{M}}$, respectively; and $K_{\mathscr{M}}$ be the bound for the sectional curvature of the manifold $\mathscr{M}$, then, we have
\begin{eqnarray}\label{eq: spectrum bounds by graph discrete paper}
\left\vert \lambda_{\bm{L}_{\mathscr{G}_{\mathscr{M}}},i}-\lambda_{\mathscr{M},i}\right\vert
&\leq& C_{n,D_{\mathscr{M}},r_{\mathscr{M}}}\left[(\epsilon/\kappa + K_{\mathscr{M}}\kappa^2)\lambda_{\mathscr{M},i}+\kappa \lambda_{\mathscr{M},i}^{3/2}\right],
\end{eqnarray}
where $\lambda_{\mathscr{M},i}$ are eigenvalues of the Laplace-Beltrami operator on the manifold $\mathscr{M}$ and $C_{n,D_{\mathscr{M}},r_{\mathscr{M}}}$ is the constant for the $i$-th eigenvalue associated to the underlying manifold properties of $\mathscr{M}$. The manifold $
mathscr{M}$ properties are its diameter $D_{\mathscr{M}}$ and its injectivity radius $r_{\mathscr{M}}$. 

Moreover, by setting the transition matrix $\bm{P}_{\mathscr{G}_{\mathscr{M}}}\define\bm{D}^{-1}_{\mathscr{G}_{\mathscr{M}}}\bm{A}_{\mathscr{G}_{\mathscr{M}}}$ associated to the approximation graph $\mathscr{G}_{\mathscr{M}}$, and let $\lambda_{\bm{P}_{\mathscr{G}_{\mathscr{M}}},i}$ be the $i$-th eigenvalue of the transition matrix $\bm{P}_{\mathscr{G}_{\mathscr{M}}}$, then, we have
\begin{eqnarray}\label{eq: lambda P g M and lambda L g m relation}
\lambda_{\bm{P}_{\mathscr{G}_{\mathscr{M}}},i}&=&1-\frac{\lambda_{\bm{L}_{\mathscr{G}_{\mathscr{M}}},i}}{\sum\limits_{j=1}^N \omega_{i,j}},
\end{eqnarray}
where $\omega_{i,j}$ is definied by Eq.~\eqref{eq: edge weight def}.
\end{theorem}

\section{Tail Bounds for Functions of Weighted Tensor Sums}\label{sec: Tail Bounds for Functions of Weighted Tensor Sums}

In this section, we will derive tail bounds for functions of weighted tensor sums derived from random walks on Riemannian manifolds. In Section~\ref{sec: exp evaluation}, the key quantity about the expectation for the norm of the weighted sum of expander tensors will be determined. In Section~\ref{sec: Upper Bound}, the upper bounds for functions of weighted tensor sums derived from random walks on Riemannian manifolds will be given. On the other hand, the lower bounds for functions of weighted tensor sums derived from random walks on Riemannian manifolds will also be given in Section~\ref{sec: Lower Bound}. 

\subsection{Evaluation $\mathbb{E}_{v_i \in \mathscr{V}}\left[\left\Vert\sum\limits_{i=1}^{\ell}w_i \Psi_i \left(g (f(v_i))\right)\right\Vert\right]$}\label{sec: exp evaluation}

Let $\mathscr{G}_{\mathscr{M}}(\epsilon,\mu,\kappa)=(\mathscr{V}, \mathscr{E}, \mathscr{W})$ be the approximation graph for $\mathscr{M}$ with transition matrix $\bm{P}_{\mathscr{G}_{\mathscr{M}}}\in \mathbb{R}^{N \times N}$, and $f: \mathscr{V} \rightarrow \in \mathbb{C}^{I_1 \times \cdots \times I_M \times I_1 \times \cdots \times I_M}$ be a Hermitian tensor-valued function. If we assume that the starting vertex is selected uniformly among all $N$ vertices from the graph $\mathscr{G}_{\mathscr{M}}$, then, we have
\begin{eqnarray}\label{eq1: sec: exp evaluation}
\lefteqn{\mathbb{E}_{v_i \in \mathscr{V}}\left[\left\Vert\sum\limits_{i=1}^{\ell}w_i \Psi_i \left(g (f(v_i))\right)\right\Vert\right]}\nonumber \\
&=&\sum\limits_{v_0,v_1,\ldots, v_\ell\in \mathscr{V}}\mathrm{P}(v_0)\mathrm{P}(v_0,v_1)\mathrm{P}(v_1,v_2)\ldots \mathrm{P}(v_{\ell-1},v_\ell)\left\Vert\sum\limits_{i=1}^{\ell}w_i \Psi_i \left(g (f(v_i))\right)\right\Vert\nonumber \\
&=&\frac{1}{N}\sum\limits_{v_0,v_1,\ldots, v_\ell\in \mathscr{V}}\mathrm{P}(v_0,v_1)\mathrm{P}(v_1,v_2)\ldots \mathrm{P}(v_{\ell-1},v_\ell)\left\Vert\sum\limits_{i=1}^{\ell}w_i \Psi_i \left(g (f(v_i))\right)\right\Vert\nonumber \\
&\leq&\frac{1}{N}\sum\limits_{v_0,v_1,\ldots, v_\ell\in \mathscr{V}}\mathrm{P}(v_0,v_1)\mathrm{P}(v_1,v_2)\ldots \mathrm{P}(v_{\ell-1},v_\ell)\sum\limits_{i=1}^{\ell}w_i\left\Vert\Psi_i \left(g (f(v_i))\right)\right\Vert\nonumber \\
&=&\frac{1}{N}\sum\limits_{j=1}^{N}\sum\limits_{i=1}^{\ell}\underbrace{\left(\sum\limits_{l=1}^{N}\left[\bm{P}_{\mathscr{G}_{\mathscr{M}}}^{(i)}\right]_{l,j}\right)}_{\define \mathscr{C}_{\bm{P}_{\mathscr{G}_{\mathscr{M}}}^{(i)},j}}w_i\left\Vert\Psi_j \left(g (f(v_j))\right)\right\Vert\nonumber \\
&=&\frac{1}{N}\sum\limits_{j=1}^{N}\sum\limits_{i=1}^{\ell}\mathscr{C}_{\bm{P}_{\mathscr{G}_{\mathscr{M}}}^{(i)},j}w_i\left\Vert\Psi_j \left(g (f(v_j))\right)\right\Vert
\end{eqnarray}
where $\mathrm{P}(v_i, v_{i+1})$ is the one step transition probability from the vertex $v_i$ to the vertex $v_{i+1}$, and $\bm{P}_{\mathscr{G}_{\mathscr{M}}}^{(i)}$ represents the $i$-steps transition matrix based on the transition matrix $\bm{P}_{\mathscr{G}_{\mathscr{M}}}$. Note that the term $\left(\sum\limits_{l=1}^{N}\left[\bm{P}_{\mathscr{G}_{\mathscr{M}}}^{(i)}\right]_{l,j}\right)\define \mathscr{C}_{\bm{P}_{\mathscr{G}_{\mathscr{M}}}^{(i)},j}$ is the column sum for the $i$-steps transition matrix. 

On the other hand, if we have $w_i>0$ and $\left\Vert\Psi_i \left(g (f(v_i))\right)\right\Vert>0$ for all $i=1,2\ldots,\ell$, we have
\begin{eqnarray}\label{eq2: sec: exp evaluation}
\lefteqn{\mathbb{E}_{v_i \in \mathscr{V}}\left[\left\Vert\sum\limits_{i=1}^{\ell}w_i \Psi_i \left(g (f(v_i))\right)\right\Vert\right]}\nonumber \\
&=&\sum\limits_{v_0,v_1,\ldots, v_\ell\in \mathscr{V}}\mathrm{P}(v_0)\mathrm{P}(v_0,v_1)\mathrm{P}(v_1,v_2)\ldots \mathrm{P}(v_{\ell-1},v_\ell)\left\Vert\sum\limits_{i=1}^{\ell}w_i \Psi_i \left(g (f(v_i))\right)\right\Vert\nonumber \\
&=&\frac{1}{N}\sum\limits_{v_0,v_1,\ldots, v_\ell\in \mathscr{V}}\mathrm{P}(v_0,v_1)\mathrm{P}(v_1,v_2)\ldots \mathrm{P}(v_{\ell-1},v_\ell)\left\Vert\sum\limits_{i=1}^{\ell}w_i \Psi_i \left(g (f(v_i))\right)\right\Vert\nonumber \\
&\geq_1&\frac{1}{N}\sum\limits_{v_0,v_1,\ldots, v_\ell\in \mathscr{V}}\mathrm{P}(v_0,v_1)\mathrm{P}(v_1,v_2)\ldots \mathrm{P}(v_{\ell-1},v_\ell)\sum\limits_{i=1}^{\ell}w_i c_i \left\Vert\Psi_i \left(g (f(v_i))\right)\right\Vert\nonumber \\
&=&\frac{1}{N}\sum\limits_{j=1}^{N}\sum\limits_{i=1}^{\ell}\underbrace{\left(\sum\limits_{l=1}^{N}\left[\bm{P}_{\mathscr{G}_{\mathscr{M}}}^{(i)}\right]_{l,j}\right)}_{\define \mathscr{C}_{\bm{P}_{\mathscr{G}_{\mathscr{M}}}^{(i)},j}}w_i c_i \left\Vert\Psi_j \left(g (f(v_j))\right)\right\Vert\nonumber \\
&=&\frac{1}{N}\sum\limits_{j=1}^{N}\sum\limits_{i=1}^{\ell}\mathscr{C}_{\bm{P}_{\mathscr{G}_{\mathscr{M}}}^{(i)},j}w_i c_i \left\Vert\Psi_j \left(g (f(v_j))\right)\right\Vert,
\end{eqnarray}
where $\geq_1$ is valid if we set $c_i$ as
\begin{eqnarray}\label{eq3: sec: exp evaluation}
c_i &\leq& \min\limits_{\forall v_j \in \mathscr{V}}\frac{\frac{1}{\ell}\left\Vert\sum\limits_{j=1}^{\ell}w_j \Psi_j \left(g (f(v_j))\right)\right\Vert}{w_i\left\Vert\Psi_i \left(g (f(v_i))\right)\right\Vert}.  
\end{eqnarray}

\subsection{Upper Bound}\label{sec: Upper Bound}

In this section, we will apply Lemma~\ref{lma: column sum UB-1} to derive Theorem~\ref{thm: Upper Bound-1} and Lemma~\ref{lma: column sum UB-2} to derive Theorem~\ref{thm: Upper Bound-2}. Both Theorem~\ref{thm: Upper Bound-1} and Theorem~\ref{thm: Upper Bound-2} provide the upper right tail bounds for functions of weighted tensor sums derived from random walks on Riemannian manifolds.

\begin{theorem}\label{thm: Upper Bound-1}
Given a manifold $\mathscr{M}$ with dimension $n$ such that the sectional curvature of the manifold $\mathscr{M}$ is bounded by $K_{\mathscr{M}}$ with the diameter $D_{\mathscr{M}}$ and the injectivity radius $r_{\mathscr{M}}$.  The approximation graph for $\mathscr{M}$ is constructed by Theorem~\ref{thm: discrete graph Laplacian appro by Laplace-Beltrami Operator on Manifold} as $\mathscr{G}_{\mathscr{M}}(\epsilon,\mu,\kappa)=(\mathscr{V}, \mathscr{E}, \mathscr{W})$ with $\left\vert\mathscr{V}\right\vert=N$. Let the funciton $f: \mathscr{V} \rightarrow \in \mathbb{C}^{I_1 \times \cdots \times I_M \times I_1 \times \cdots \times I_M}$ be a Hermitian tensor-valued function such that all Hermitian tensors under the function $f$ have eigenvalues within in the real interval $[c, d]$. We assume that $\max\limits_{j \in \{1,2,\ldots N\}}\mathscr{C}_{\bm{P}_{\mathscr{G}_{\mathscr{M}}}^{(i)},j}\geq N-1$, where $i=1,2,\ldots,\ell$. $\Psi$ is a normalized positive linear map and a probability vector with size $\ell$ is given. For continuous real functions $g,h$ with $g$ convexity in the real interval $[c, d]$ and $h(s)>0, m_g s + b_{g,U}>0, m_g s + b_{g,L}>0$ for $s \in [c, d]$, where $m_g, b_{g,U}$ and $b_{g,L}$ are provided by Eq.~\eqref{eq: m g b g U b g L def}. Consider a random walk with $\ell$ steps over $\mathscr{M}$, then, for any $\theta>0$, we have
\begin{eqnarray}\label{eq1: thm: Upper Bound-1}
\mathrm{P}\left( \left\Vert h\left(\sum\limits_{i=1}^\ell w_i \Psi(f(v_i))\right) \right\Vert \geq \theta\right)\leq \frac{\sum\limits_{j=1}^{N}\sum\limits_{i=1}^{\ell}\left(N -\left\vert\left(\frac{C_1}{\sum\limits_{k=1}^N \omega_{2,k}}\right)^i\right\vert\right)w_i\left\Vert\Psi \left(g (f(v_j))\right)\right\Vert}{N\theta\min\limits_{c \leq s \leq d}\left[\frac{m_g s + b_{g,L}}{h(s)}\right]},
\end{eqnarray}
where $C_1$ is defined by
\begin{eqnarray}
C_1&\define&\min\{\lambda_{\mathscr{M},2}+C_{n,D_{\mathscr{M}},r_{\mathscr{M}}}\left[(\epsilon/\kappa + K_{\mathscr{M}}\kappa^2)\lambda_{\mathscr{M},2}+\kappa \lambda_{\mathscr{M},2}^{3/2}\right],\nonumber \\
&&\lambda_{\mathscr{M},2}-C_{n,D_{\mathscr{M}},r_{\mathscr{M}}}\left[(\epsilon/\kappa + K_{\mathscr{M}}\kappa^2)\lambda_{\mathscr{M},2}+\kappa \lambda_{\mathscr{M},2}^{3/2}\right]\}.
\end{eqnarray}
\end{theorem}
\textbf{Proof:}
Since we have
\begin{eqnarray}\label{eq2: thm: Upper Bound-1}
\mathrm{P}\left( \left\Vert h\left(\sum\limits_{i=1}^\ell w_i \Psi(f(v_i))\right) \right\Vert \geq \theta\right)&\leq_1&\frac{\mathbb{E}\left[\left\Vert h\left(\sum\limits_{i=1}^\ell w_i \Psi(f(v_i))\right) \right\Vert\right]}{\theta}\nonumber \\
&\leq_2& \frac{\mathbb{E}\left[\left\Vert \sum\limits_{i=1}^\ell w_i \Psi(g(f(v_i))) \right\Vert\right]}{\theta\min\limits_{c \leq s \leq d}\left[\frac{m_g s + b_{g,L}}{h(s)}\right]},
\end{eqnarray}
where $\leq_1$ comes from Markov inequality, and $\leq_2$ comes from Theorem~\ref{thm: Tensor Inequalities by MP}. From  Eq.~\eqref{eq1: sec: exp evaluation} and Eq.~\eqref{eq2: thm: Upper Bound-1}, we have
\begin{eqnarray}\label{eq3: thm: Upper Bound-1}
\mathrm{P}\left( \left\Vert h\left(\sum\limits_{i=1}^\ell w_i \Psi(f(v_i))\right) \right\Vert \geq \theta\right)&\leq&\frac{\sum\limits_{j=1}^{N}\sum\limits_{i=1}^{\ell}\mathscr{C}_{\bm{P}_{\mathscr{G}_{\mathscr{M}}}^{(i)},j}w_i\left\Vert\Psi \left(g (f(v_j))\right)\right\Vert}{N\theta\min\limits_{c \leq s \leq d}\left[\frac{m_g s + b_{g,L}}{h(s)}\right]} \nonumber \\
&\leq_1& \frac{\sum\limits_{j=1}^{N}\sum\limits_{i=1}^{\ell}\left(N - \left\vert\lambda^i_2 (\bm{P}_{\mathscr{G}_{\mathscr{M}}})\right\vert\right)w_i\left\Vert\Psi \left(g (f(v_j))\right)\right\Vert}{N\theta\min\limits_{c \leq s \leq d}\left[\frac{m_g s + b_{g,L}}{h(s)}\right]},
\end{eqnarray}
where $\leq_1$ comes from Lemma~\ref{lma: column sum UB-1} and $\lambda_2(\bm{P}_{\mathscr{G}_{\mathscr{M}}})$ is the second largest absolute eigenvaule of the transition matrix $\bm{P}_{\mathscr{G}_{\mathscr{M}}}$. 

Then, we have the desired result by applying Eq.~\eqref{eq: spectrum bounds by graph discrete paper} to $\lambda_2 (\bm{P}_{\mathscr{G}_{\mathscr{M}}})$.
$\hfill\Box$

The next Theorem~\ref{thm: Upper Bound-2} will assume the different condition for the transition matrix $\bm{P}_{\mathscr{G}_{\mathscr{M}}}$.

\begin{theorem}\label{thm: Upper Bound-2}
Given a manifold $\mathscr{M}$ with dimension $n$ such that the sectional curvature of the manifold $\mathscr{M}$ is bounded by $K_{\mathscr{M}}$ with the diameter $D_{\mathscr{M}}$ and the injectivity radius $r_{\mathscr{M}}$.  The approximation graph for $\mathscr{M}$ is constructed by Theorem~\ref{thm: discrete graph Laplacian appro by Laplace-Beltrami Operator on Manifold} as $\mathscr{G}_{\mathscr{M}}(\epsilon,\mu,\kappa)=(\mathscr{V}, \mathscr{E}, \mathscr{W})$ with $\left\vert\mathscr{V}\right\vert=N$. Let the funciton $f: \mathscr{V} \rightarrow \in \mathbb{C}^{I_1 \times \cdots \times I_M \times I_1 \times \cdots \times I_M}$ be a Hermitian tensor-valued function such that all Hermitian tensors under the function $f$ have eigenvalues within in the real interval $[c, d]$. Twice of any diagonal entry of $\bm{P}_{\mathscr{G}_{\mathscr{M}}}$ is assumed to be greater than one. $\Psi$ is a normalized positive linear map and a probability vector with size $\ell$ is given. For continuous real functions $g,h$ with $g$ convexity in the real interval $[c, d]$ and $h(s)>0, m_g s + b_{g,U}>0, m_g s + b_{g,L}>0$ for $s \in [c, d]$, where $m_g, b_{g,U}$ and $b_{g,L}$ are provided by Eq.~\eqref{eq: m g b g U b g L def}. Consider a random walk with $\ell$ steps over $\mathscr{M}$, then, for any $\theta>0$ and $q \in (0,1)$, we have
\begin{eqnarray}\label{eq1: thm: Upper Bound-2}
\mathrm{P}\left( \left\Vert h\left(\sum\limits_{i=1}^\ell w_i \Psi(f(v_i))\right) \right\Vert \geq \theta\right)&\leq&\frac{\sum\limits_{j=1}^{N}\sum\limits_{i=1}^{\ell}\left[\frac{\min\limits_{k \in \{1,2,\ldots,N\}}\left\vert C^i_{2,k}\right\vert}{\left(2\left[\bm{P}^i_{\mathscr{G}_{\mathscr{M}}}\right]_{j,j} - 1\right)^q}\right]^{\frac{1}{1-q}}w_i\left\Vert\Psi \left(g (f(v_j))\right)\right\Vert}{N\theta\min\limits_{c \leq s \leq d}\left[\frac{m_g s + b_{g,L}}{h(s)}\right]}
\end{eqnarray}
where $\left[\bm{P}^i_{\mathscr{G}_{\mathscr{M}}}\right]_{j,j}$ is the $(j,j)$ entry in the transition matrix $\bm{P}^i_{\mathscr{G}_{\mathscr{M}}}$, and $C_{2,k}$ is defined by
\begin{eqnarray}
C_{2,k}&\define&\max\Bigg\{1-\frac{\lambda_{\mathscr{M},k}+C_{n,D_{\mathscr{M}},r_{\mathscr{M}}}\left[(\epsilon/\kappa + K_{\mathscr{M}}\kappa^2)\lambda_{\mathscr{M},k}+\kappa \lambda_{\mathscr{M},k}^{3/2}\right]}{\sum\limits_{j=1}^N \omega_{k,j}},\nonumber \\
&&1-\frac{\lambda_{\mathscr{M},k}-C_{n,D_{\mathscr{M}},r_{\mathscr{M}}}\left[(\epsilon/\kappa + K_{\mathscr{M}}\kappa^2)\lambda_{\mathscr{M},k}+\kappa \lambda_{\mathscr{M},k}^{3/2}\right]}{\sum\limits_{j=1}^N \omega_{k,j}}\Bigg\}.
\end{eqnarray}
\end{theorem}
\textbf{Proof:}
Because we have
\begin{eqnarray}\label{eq2: thm: Upper Bound-2}
\mathrm{P}\left( \left\Vert h\left(\sum\limits_{i=1}^\ell w_i \Psi(f(v_i))\right) \right\Vert \geq \theta\right)&\leq_1&\frac{\mathbb{E}\left[\left\Vert h\left(\sum\limits_{i=1}^\ell w_i \Psi(f(v_i))\right) \right\Vert\right]}{\theta}\nonumber \\
&\leq_2& \frac{\mathbb{E}\left[\left\Vert \sum\limits_{i=1}^\ell w_i \Psi(g(f(v_i))) \right\Vert\right]}{\theta\min\limits_{c \leq s \leq d}\left[\frac{m_g s + b_{g,L}}{h(s)}\right]},
\end{eqnarray}
where $\leq_1$ comes from Markov inequality, and $\leq_2$ comes from Theorem~\ref{thm: Tensor Inequalities by MP}. From  Eq.~\eqref{eq1: sec: exp evaluation} and Eq.~\eqref{eq2: thm: Upper Bound-1}, we have
\begin{eqnarray}\label{eq3: thm: Upper Bound-2}
\mathrm{P}\left( \left\Vert h\left(\sum\limits_{i=1}^\ell w_i \Psi(f(v_i))\right) \right\Vert \geq \theta\right)&\leq&\frac{\sum\limits_{j=1}^{N}\sum\limits_{i=1}^{\ell}\mathscr{C}_{\bm{P}_{\mathscr{G}_{\mathscr{M}}}^{(i)},j}w_i\left\Vert\Psi \left(g (f(v_j))\right)\right\Vert}{N\theta\min\limits_{c \leq s \leq d}\left[\frac{m_g s + b_{g,L}}{h(s)}\right]} \nonumber \\
&\leq_1& \frac{\sum\limits_{j=1}^{N}\sum\limits_{i=1}^{\ell}\left[\frac{\min\limits_{k \in \{1,2,\ldots,N\}}\left\vert\lambda^i_{k}\left(\bm{P}_{\mathscr{G}_{\mathscr{M}}}\right)\right\vert}{\left(2\left[\bm{P}^i_{\mathscr{G}_{\mathscr{M}}}\right]_{j,j} - 1\right)^q}\right]^{\frac{1}{1-q}}w_i\left\Vert\Psi \left(g (f(v_j))\right)\right\Vert}{N\theta\min\limits_{c \leq s \leq d}\left[\frac{m_g s + b_{g,L}}{h(s)}\right]},
\end{eqnarray}
where $\leq_1$ comes from Lemma~\ref{lma: column sum UB-2}.

Then, we have the desired result by applying Eq.~\eqref{eq: spectrum bounds by graph discrete paper} and Theorem~\ref{thm: discrete graph Laplacian appro by Laplace-Beltrami Operator on Manifold} to $\lambda_k(\bm{P}_{\mathscr{G}_{\mathscr{M}}})$.
$\hfill\Box$

\subsection{Lower Bound}\label{sec: Lower Bound}

In this section, we will apply Lemma~\ref{lma: column sum LB-1} to derive Theorem~\ref{thm: Lower Bound-1} and Lemma~\ref{lma: column sum LB-2} to derive Theorem~\ref{thm: Lower Bound-2}. Both Theorem~\ref{thm: Lower Bound-1} and Theorem~\ref{thm: Lower Bound-2} provide the lower right tail bounds for functions of weighted tensor sums derived from random walks on Riemannian manifolds.

\begin{theorem}\label{thm: Lower Bound-1}
Given a manifold $\mathscr{M}$ with dimension $n$ such that the sectional curvature of the manifold $\mathscr{M}$ is bounded by $K_{\mathscr{M}}$ with the diameter $D_{\mathscr{M}}$ and the injectivity radius $r_{\mathscr{M}}$.  The approximation graph for $\mathscr{M}$ is constructed by Theorem~\ref{thm: discrete graph Laplacian appro by Laplace-Beltrami Operator on Manifold} as $\mathscr{G}_{\mathscr{M}}(\epsilon,\mu,\kappa)=(\mathscr{V}, \mathscr{E}, \mathscr{W})$ with $\left\vert\mathscr{V}\right\vert=N$. Let the funciton $f: \mathscr{V} \rightarrow \in \mathbb{C}^{I_1 \times \cdots \times I_M \times I_1 \times \cdots \times I_M}$ be a Hermitian tensor-valued function such that all Hermitian tensors under the function $f$ have eigenvalues within in the real interval $[c, d]$.  All diagonal entries of $\bm{P}_{\mathscr{G}_{\mathscr{M}}}$ are assumed less than one. $\Psi$ is a normalized positive linear map and a probability vector with size $\ell$ with positive entries is given. For continuous real functions $g,h$ with $g$ convexity in the real interval $[c, d]$ and $h(s)>0, m_g s + b_{g,U}>0, m_g s + b_{g,L}>0$ for $s \in [c, d]$, where $m_g, b_{g,U}$ and $b_{g,L}$ are provided by Eq.~\eqref{eq: m g b g U b g L def}. We assume that the norm for any Hermitian $\mathcal{X}$ satisfies $\left\Vert h(\mathcal{X}) \right\Vert \leq \mathrm{A}$ for some positive number $\mathrm{A}$. Consider a random walk with $\ell$ steps over $\mathscr{M}$, then, for any $\theta>0$ and $q \in (0,1)$, we have
\begin{eqnarray}\label{eq1: thm: Lower Bound-1}
\lefteqn{\mathrm{P}\left( \left\Vert h\left(\sum\limits_{i=1}^\ell w_i \Psi(f(v_i))\right) \right\Vert \geq \theta\right)}\nonumber \\
&\geq&\frac{\sum\limits_{j=1}^{N}\sum\limits_{i=1}^{\ell}\left[ [\bm{P}^i_{\mathscr{G}_{\mathscr{M}}}]_{j,j}+\frac{C_{3,i,j}^{1/(1-q)}}{(1-[\bm{P}^i_{\mathscr{G}_{\mathscr{M}}}]_{j,j})^{q/(1-q)}}\right]w_i c_i \left\Vert\Psi \left(g (f(v_j))\right)\right\Vert - \theta\max\limits_{c \leq s \leq d}\left[\frac{m_g s + b_{g,U}}{h(s)}\right]}{N(\mathrm{A}-\theta)\max\limits_{c \leq s \leq d}\left[\frac{m_g s + b_{g,U}}{h(s)}\right]},
\end{eqnarray}
where $C_{3,i,j}$ is defined by
\begin{eqnarray}\label{eq1.5: thm: Lower Bound-1}
C_{3,i,j}&\define&\min \Bigg\{\left\vert\left\{1-\frac{\lambda_{\mathscr{M},j}+C_{n,D_{\mathscr{M}},r_{\mathscr{M}}}\left[(\epsilon/\kappa + K_{\mathscr{M}}\kappa^2)\lambda_{\mathscr{M},j}+\kappa \lambda_{\mathscr{M},j}^{3/2}\right]}{\sum\limits_{k=1}^N \omega_{j,k}}\right\}^i - [\bm{P}^i_{\mathscr{G}_{\mathscr{M}}}]_{j,j}\right\vert,\nonumber \\
&&\left\vert\left\{1-\frac{\lambda_{\mathscr{M},j}-C_{n,D_{\mathscr{M}},r_{\mathscr{M}}}\left[(\epsilon/\kappa + K_{\mathscr{M}}\kappa^2)\lambda_{\mathscr{M},j}+\kappa \lambda_{\mathscr{M},j}^{3/2}\right]}{\sum\limits_{k=1}^N \omega_{j,k}}\right\}^i - [\bm{P}^i_{\mathscr{G}_{\mathscr{M}}}]_{j,j}\right\vert\Bigg\}.
\end{eqnarray}
\end{theorem}
\textbf{Proof:}
Because we have
\begin{eqnarray}\label{eq2: thm: Lower Bound-1}
\mathrm{P}\left( \left\Vert h\left(\sum\limits_{i=1}^\ell w_i \Psi(f(v_i))\right) \right\Vert \geq \theta\right)&\geq_1&\frac{\mathbb{E}\left[\left\Vert h\left(\sum\limits_{i=1}^\ell w_i \Psi(f(v_i))\right) \right\Vert\right]- \theta}{\mathrm{A} - \theta}\nonumber \\
&\geq_2& \frac{\mathbb{E}\left[\left\Vert \sum\limits_{i=1}^\ell w_i \Psi(g(f(v_i))) \right\Vert\right] - \theta\max\limits_{c \leq s \leq d}\left[\frac{m_g s + b_{g,U}}{h(s)}\right] }{(\mathrm{A}-\theta)\max\limits_{c \leq s \leq d}\left[\frac{m_g s + b_{g,U}}{h(s)}\right]},
\end{eqnarray}
where $\geq_1$ comes from reverse Markov inequality, and $\geq_2$ comes from Theorem~\ref{thm: Tensor Inequalities by MP}. From  Eq.~\eqref{eq2: sec: exp evaluation} and Eq.~\eqref{eq2: thm: Lower Bound-1}, we have
\begin{eqnarray}\label{eq3: thm: Lower Bound-1}
\mathrm{P}\left( \left\Vert h\left(\sum\limits_{i=1}^\ell w_i \Psi(f(v_i))\right) \right\Vert \geq \theta\right)\geq\frac{\sum\limits_{j=1}^{N}\sum\limits_{i=1}^{\ell}\mathscr{C}_{\bm{P}_{\mathscr{G}_{\mathscr{M}}}^{(i)},j}w_i c_i \left\Vert\Psi \left(g (f(v_j))\right)\right\Vert - \theta\max\limits_{c \leq s \leq d}\left[\frac{m_g s + b_{g,U}}{h(s)}\right]}{N(\mathrm{A}-\theta)\max\limits_{c \leq s \leq d}\left[\frac{m_g s + b_{g,U}}{h(s)}\right]} \nonumber \\
\geq_1\frac{\sum\limits_{j=1}^{N}\sum\limits_{i=1}^{\ell}\left[ [\bm{P}^i_{\mathscr{G}_{\mathscr{M}}}]_{j,j}+\frac{\left\vert\lambda_j(\bm{P}^i_{\mathscr{G}_{\mathscr{M}}}) - [\bm{P}^i_{\mathscr{G}_{\mathscr{M}}}]_{j,j}\right\vert^{1/(1-q)}}{(1-[\bm{P}^i_{\mathscr{G}_{\mathscr{M}}}]_{j,j})^{q/(1-q)}}\right]w_i c_i \left\Vert\Psi \left(g (f(v_j))\right)\right\Vert - \theta\max\limits_{c \leq s \leq d}\left[\frac{m_g s + b_{g,U}}{h(s)}\right]}{N(\mathrm{A}-\theta)\max\limits_{c \leq s \leq d}\left[\frac{m_g s + b_{g,U}}{h(s)}\right]},
\end{eqnarray}
where $\geq_1$ comes from Lemma~\ref{lma: column sum LB-1}.

Then, we have the desired result by applying Eq.~\eqref{eq: spectrum bounds by graph discrete paper} and Theorem~\ref{thm: discrete graph Laplacian appro by Laplace-Beltrami Operator on Manifold} to $\lambda_j(\bm{P}^i_{\mathscr{G}_{\mathscr{M}}})$.
$\hfill\Box$

\begin{theorem}\label{thm: Lower Bound-2}
Given a manifold $\mathscr{M}$ with dimension $n$ such that the sectional curvature of the manifold $\mathscr{M}$ is bounded by $K_{\mathscr{M}}$ with the diameter $D_{\mathscr{M}}$ and the injectivity radius $r_{\mathscr{M}}$.  The approximation graph for $\mathscr{M}$ is constructed by Theorem~\ref{thm: discrete graph Laplacian appro by Laplace-Beltrami Operator on Manifold} as $\mathscr{G}_{\mathscr{M}}(\epsilon,\mu,\kappa)=(\mathscr{V}, \mathscr{E}, \mathscr{W})$ with $\left\vert\mathscr{V}\right\vert=N$. Let the funciton $f: \mathscr{V} \rightarrow \in \mathbb{C}^{I_1 \times \cdots \times I_M \times I_1 \times \cdots \times I_M}$ be a Hermitian tensor-valued function such that all Hermitian tensors under the function $f$ have eigenvalues within in the real interval $[c, d]$. All diagonal entries of $\bm{P}_{\mathscr{G}_{\mathscr{M}}}$ are assumed less than one. $\Psi$ is a normalized positive linear map and a probability vector with positive entries and size $\ell$ is given. For continuous real functions $g,h$ with $g$ convexity in the real interval $[c, d]$ and $h(s)>0, m_g s + b_{g,U}>0, m_g s + b_{g,L}>0$ for $s \in [c, d]$, where $m_g, b_{g,U}$ and $b_{g,L}$ are provided by Eq.~\eqref{eq: m g b g U b g L def}. We assume that the norm for any Hermitian $\mathcal{X}$ satisfies $\left\Vert h(\mathcal{X}) \right\Vert \leq \mathrm{A}$ for some positive number $\mathrm{A}$. Consider a random walk with $\ell$ steps over $\mathscr{M}$, then, for any $\theta>0$ and $q \in (0,1)$, we have
\begin{eqnarray}\label{eq1: thm: Lower Bound-2}
\lefteqn{\mathrm{P}\left( \left\Vert h\left(\sum\limits_{i=1}^\ell w_i \Psi(f(v_i))\right) \right\Vert \geq \theta\right)}\nonumber \\
&\geq&\frac{\sum\limits_{j=1}^{N}\sum\limits_{i=1}^{\ell}C_{4,i,j}^2w_i c_i \left\Vert\Psi \left(g (f(v_j))\right)\right\Vert - \theta\max\limits_{c \leq s \leq d}\left[\frac{m_g s + b_{g,U}}{h(s)}\right]}{N(\mathrm{A}-\theta)\max\limits_{c \leq s \leq d}\left[\frac{m_g s + b_{g,U}}{h(s)}\right]},
\end{eqnarray}
where $C_{4,i,j}$ is defined by
\begin{eqnarray}\label{eq1.5: thm: Lower Bound-2}
C_{4,i,j}&\define&\min \Bigg\{\left\vert\left\{1-\frac{\lambda_{\mathscr{M},j}+C_{n,D_{\mathscr{M}},r_{\mathscr{M}}}\left[(\epsilon/\kappa + K_{\mathscr{M}}\kappa^2)\lambda_{\mathscr{M},j}+\kappa \lambda_{\mathscr{M},j}^{3/2}\right]}{\sum\limits_{k=1}^N \omega_{j,k}}\right\}^i\right\vert,\nonumber \\
&&\left\vert\left\{1-\frac{\lambda_{\mathscr{M},j}-C_{n,D_{\mathscr{M}},r_{\mathscr{M}}}\left[(\epsilon/\kappa + K_{\mathscr{M}}\kappa^2)\lambda_{\mathscr{M},j}+\kappa \lambda_{\mathscr{M},j}^{3/2}\right]}{\sum\limits_{k=1}^N \omega_{j,k}}\right\}^i\right\vert\Bigg\}.
\end{eqnarray}
\end{theorem}
\textbf{Proof:}
Because we have
\begin{eqnarray}\label{eq2: thm: Lower Bound-2}
\mathrm{P}\left( \left\Vert h\left(\sum\limits_{i=1}^\ell w_i \Psi(f(v_i))\right) \right\Vert \geq \theta\right)&\geq_1&\frac{\mathbb{E}\left[\left\Vert h\left(\sum\limits_{i=1}^\ell w_i \Psi(f(v_i))\right) \right\Vert\right]- \theta}{\mathrm{A} - \theta}\nonumber \\
&\geq_2& \frac{\mathbb{E}\left[\left\Vert \sum\limits_{i=1}^\ell w_i \Psi(g(f(v_i))) \right\Vert\right] - \theta\max\limits_{c \leq s \leq d}\left[\frac{m_g s + b_{g,U}}{h(s)}\right] }{(\mathrm{A}-\theta)\max\limits_{c \leq s \leq d}\left[\frac{m_g s + b_{g,U}}{h(s)}\right]},
\end{eqnarray}
where $\geq_1$ comes from reverse Markov inequality, and $\geq_2$ comes from Theorem~\ref{thm: Tensor Inequalities by MP}. From  Eq.~\eqref{eq2: sec: exp evaluation} and Eq.~\eqref{eq2: thm: Lower Bound-2}, we have
\begin{eqnarray}\label{eq3: thm: Lower Bound-2}
\mathrm{P}\left( \left\Vert h\left(\sum\limits_{i=1}^\ell w_i \Psi(f(v_i))\right) \right\Vert \geq \theta\right)\geq\frac{\sum\limits_{j=1}^{N}\sum\limits_{i=1}^{\ell}\mathscr{C}_{\bm{P}_{\mathscr{G}_{\mathscr{M}}}^{(i)},j}w_i c_i \left\Vert\Psi \left(g (f(v_j))\right)\right\Vert - \theta\max\limits_{c \leq s \leq d}\left[\frac{m_g s + b_{g,U}}{h(s)}\right]}{N(\mathrm{A}-\theta)\max\limits_{c \leq s \leq d}\left[\frac{m_g s + b_{g,U}}{h(s)}\right]} \nonumber \\
\geq_1\frac{\sum\limits_{j=1}^{N}\sum\limits_{i=1}^{\ell}\left\vert\lambda_j(\bm{P}^i_{\mathscr{G}_{\mathscr{M}}})\right\vert^2w_i c_i \left\Vert\Psi \left(g (f(v_j))\right)\right\Vert - \theta\max\limits_{c \leq s \leq d}\left[\frac{m_g s + b_{g,U}}{h(s)}\right]}{N(\mathrm{A}-\theta)\max\limits_{c \leq s \leq d}\left[\frac{m_g s + b_{g,U}}{h(s)}\right]},
\end{eqnarray}
where $\geq_1$ comes from Lemma~\ref{lma: column sum LB-2}.

Then, we have the desired result by applying Eq.~\eqref{eq: spectrum bounds by graph discrete paper} and Theorem~\ref{thm: discrete graph Laplacian appro by Laplace-Beltrami Operator on Manifold} to $\lambda_j(\bm{P}^i_{\mathscr{G}_{\mathscr{M}}})$.
$\hfill\Box$

\begin{remark}\label{rmk: T-product extension}
Given that tensor norms for T-product tensors can also be defined as established in~\cite{chang2021T-TensorExp}, the methodology utilized in this study can be extended to explore tail bounds for functions of weighted T-product tensor sums originating from random walks on Riemannian manifolds.
\end{remark}

\bibliographystyle{IEEETran}
\bibliography{TenExpOverManifold_Bib}

\end{document}